\documentclass[final,3p, times, 11pt]{elsarticle}

\usepackage{amssymb,amsmath,color}
\usepackage[mathscr]{eucal}
 \journal{Discrete Applied Mathematics
 }
\newtheorem{thm}{Theorem}[section]

\newtheorem{lemma}[thm]{Lemma}

\newtheorem{defin}{Definition}[section]
\newproof{pf}{Proof}

\newcommand{\bF}{ \mathbb{F}}

\begin{document}

\begin{frontmatter}



\title{Pooling designs with surprisingly high degree of error correction in a finite vector space}

\author[JG]{Jun Guo}\ead{guojun$_-$lf@163.com}
\author[KW]{Kaishun Wang\corref{cor}}
\ead{wangks@bnu.edu.cn}
\cortext[cor]{Corresponding author}

\address[JG]{Math. and Inf. College, Langfang Teachers'
College, Langfang  065000,  China }
\address[KW]{Sch. Math. Sci. \& Lab. Math. Com. Sys.,
Beijing Normal University, Beijing  100875, China}

\begin{abstract} Pooling designs are standard experimental tools in
many biotechnical applications.  It is well-known that all famous
pooling designs are constructed from mathematical structures by the
``containment matrix" method. In particular, Macula's designs (resp.
Ngo and Du's designs) are constructed by the containment relation of
subsets (resp. subspaces) in a finite set (resp. vector space).
Recently, we generalized Macula's designs and obtained a family of
pooling designs with more high degree of error correction by subsets
in a finite set. In this paper, as a generalization of Ngo and Du's
designs, we study the corresponding problems in a finite vector
space  and obtain a family of pooling designs with surprisingly high
degree of error correction. Our designs and Ngo and Du's designs
have the same number of items and pools, respectively, but the
error-tolerant property is much better than that of Ngo and Du's
designs, which was given by  D'yachkov et al. \cite{DF}, when the
dimension of the space is large enough.
\end{abstract}

\end{frontmatter}

\section{Introduction}

A group test is applicable to an arbitrary subset of clones with two
possible outcomes: a negative outcome indicates all clones in the
subset are negative, and a positive outcome indicates otherwise. A
pooling design is a specification of all tests so that they can be
performed simultaneously with the goal being to identify all
positive clones with a small number of tests \cite{CD,CD2,DH,HW}. A
pooling design is usually represented by a binary matrix with
columns indexed with items and rows indexed with pools. A cell
$(i,j)$ contains a 1-entry if and only if the $i$th pool contains
the $j$th item. By treating a column as a set of row indices
intersecting the column with a 1-entry, we can talk about the union
of several columns.
 A binary matrix is
$s^e$-{\it disjunct} if every column has at least $e+1$ 1-entries
not contained in the union of any other $s$ columns \cite{Macula2}.
An $s^0$-disjunct matrix is also called $s$-disjunct. An
$s^e$-disjunct matrix is called {\it fully $s^e$-disjunct} if it is
not $s_1^{e_1}$-disjunct whenever $s_1>s$ or $e_1>e$. An
$s^e$-disjunct matrix is $\lfloor e/2\rfloor$-error-correcting
\cite{DF}.

For positive integers $k\leq n$, let $[n]=\{1,2,\ldots,n\}$ and
$\left([n]\atop k\right)$ denote the set of all $k$-subsets of
$[n]$.

Macula \cite{Macula,Macula2} proposed a novel way of constructing
disjunct matrices by the containment relation of subsets in $[n]$.

\begin{defin}{\rm(\cite{Macula})}
For positive integers $1\leq d<k< n$, let $M(d,k,n)$ be the binary
matrix with rows indexed with $\left([n]\atop
d\right)$ and columns indexed with $\left([n]\atop k\right)$ such that $M(A,B)=1$ if
and only if $A\subseteq B$.
\end{defin}

D'yachkov et al. \cite{DF2} discussed the error-correcting property
of $M(d,k,n)$.

\begin{thm}{\rm (\cite{DF2})}\label{thm1.1}
For positive integers
$1\leq d<k< n$ and $1\leq s\leq d$, $M(d,k,n)$ is fully
$s^{e_1}$-disjunct,   where $e_1=\left(k-s\atop d-s\right)-1$.
\end{thm}

In \cite{GW}, we generalized Macula's construction and obtained a
family of pooling designs with a higher degree of error correction.

\begin{defin}{\rm(\cite{GW})}
For positive integers $1\leq d< k<n $ and $0\leq i\leq
 d$.  Let $M(i;d,k,n)$ be the binary matrix with rows indexed
with $\left([n]\atop d\right)$ and columns indexed with
$\left([n]\atop k\right)$ such that $M(A,B)=1$ if and only if
$|A\cap B|=i$.
\end{defin}

\begin{thm}\label{thm1.2}{\rm(\cite{GW})}
Let $1\leq s\leq i,\lfloor(d+1)/2\rfloor\leq i\leq d<k$ and
$n-k-s(k+d-2i)\geq d-i$. Then

\begin{itemize}

\item[\rm(i)] $M(i;d,k,n)$ is an $s^{e_2}$-disjunct matrix, where
$e_2=\left(k-s\atop i-s\right)\left(n-k-s(k+d-2i)\atop
d-i\right)-1$;

\item[\rm(ii)] For a given $k$, if $i<d,$ then $\lim\limits_{n\longrightarrow +\infty}
\frac{e_2+1}{e_1+1}=+\infty.$
\end{itemize}
\end{thm}

Now we introduce the $q$-analogue of Theorems~\ref{thm1.1} and
\ref{thm1.2}.

Let $\mathbb{F}_q$ be a finite field with $q$ elements, where $q$ is
a prime power. For a positive integer $n$, let $\mathbb{F}_q^n$ be
an $n$-dimensional vector space over $\mathbb{F}_q$. For positive
integers $k\leq n$, let $\left[[n]\atop k\right]_q$ be the set of
all $k$-dimensional subspaces of $\bF_q^{n}$. A matrix
representation of a subspace $P$ is a matrix whose rows form a basis
for $P$. When there is no danger of confusion, we use the same
symbol to denote a subspace  and its matrix representation.

Let $m_1, m_2$ be two  integers. For
brevity we use  the {\it Gaussian coefficient}
$$
\left[m_2\atop m_1\right]_q=\frac{\prod\limits_{t=m_2-m_1+1}^{m_2}(q^t-1)}{\prod\limits_{t=1}^{m_1}(q^t-1)}.
$$
By convenience   $\left[m_2\atop 0\right]_q=1$ and $\left[m_2\atop
m_1\right]_q=0$ whenever $m_1<0$  or $m_2<m_1$.
Then, by \cite{wanbook},
$$\left|\left[[n]\atop k\right]_q\right|=\left[n\atop k\right]_q.$$

Ngo and Du \cite{Ngo} constructed a family of disjunct matrices
by the containment relation of subspaces in $\bF_q^{n}$.

\begin{defin}{\rm(\cite{Ngo})}\,
For positive integers $1\leq d<k< n$, let $M_q(d,k,n)$ be the binary
matrix with rows indexed with $\left[[n]\atop d\right]_q$ and columns indexed
 with $\left[[n]\atop k\right]_q$ such that $M_q(A,B)=1$
if and only if $A\subseteq B$.
\end{defin}

D'yachkov et al. \cite{DF} discussed the error-tolerant property of $M_q(d,k,n)$.

\begin{thm}{\rm(\cite{DF})}\label{thm1.3}
For positive integers
$1\leq d<k< n, k-d\geq2$ and $1\leq \bar{s}\leq q(q^{k-1}-1)/(q^{k-d}-1)$, $M_q(d,k,n)$ is
$\bar{s}^{\bar{e}_1}$-disjunct,   where
$\bar{e}_1=q^{k-d}\left[k-1\atop d-1\right]_q-(\bar{s}-1)q^{k-d-1}\left[k-2\atop d-1\right]_q-1$.
In particular, if $\bar{s}\leq q+1$, then $M_q(d,k,n)$ is
fully $\bar{s}^{\bar{e}_1}$-disjunct.
\end{thm}

Nan and Guo \cite{NG} generalized Ngo and Du's construction and
obtained a family of   pooling designs.

\begin{defin}{\rm(\cite{NG})}
For positive integers $1\leq d< k<n $ and $\max\{0,d+k-n\}\leq i\leq
 d$.  Let $M_q(i;d,k,n)$ be the binary matrix with rows indexed with $\left[[n]\atop d\right]_q$
and columns indexed with $\left[[n]\atop k\right]_q$ such that $M_q(A,B)=1$ if and only if
$\dim(A\cap B)=i$.
\end{defin}

Note that $M_q(i;d,k,n)$ and $M_q(d,k,n)$ have the same size. In
\cite{NG}, the error-tolerant property of $M_q(i;d,k,n)$ is not well
expressed. In this paper, we discuss again the error-tolerant
property of $M_q(i;d,k,n)$.

\section{Main results}
In this section, we discuss the error-tolerant property of $M_q(i;d,k,n)$.
We begin with a useful lemma.

\begin{lemma}\label{lem2.2}
For $\max\{0,r+m-n\}\leq j\leq r$ and $ m\leq n$, let $P_0$ be a
given $m$-dimensional subspace of $\mathbb{F}_q^{n}$ and let $Q_0 $
be a given $j$-dimensional subspace of $\mathbb{F}_q^{n}$ with
$Q_0\subseteq P_0$. Then the number of $r$-dimensional subspaces of
$\mathbb{F}_q^{n}$  intersecting $P_0$ at $Q_0$ is
$f(j,r,n;m)=q^{(r-j)(m-j)}\left[n-m\atop r-j\right]_q.$ Moreover,
for the integer $0\leq\alpha\leq n+j-m-r$, the function
$f(j,r,n;m+\alpha)$ about $\alpha$ is decreasing.
\end{lemma}

\begin{pf}
Since the general linear group $GL_n(\mathbb{F}_q)$ acts
transitively on the set of such pairs $(P_0, Q_0)$, we may assume
that $P_0=(I^{(m)}\;0^{(m,n-m)}),\;Q_0=(I^{(j)}\;0^{(j,n-j)}).$ Let
$Q$ be an $r$-dimensional subspace of $\mathbb{F}_q^{n}$ satisfying
$P_0\cap Q=Q_0$. Then $Q$ has a matrix representation of the form
$$
\left(\begin{array}{ccc}
I^{(j)}&0^{(j,m-j)}&0^{(j,n-m)}\\
0^{(r-j,j)}   &A_2 & A_3
\end{array}\right),
$$
where $A_2$ is an $(r-j)\times(m-j)$ matrix and $A_3$ is an
$(r-j)$-dimensional subspace of $\mathbb{F}_q^{n-m}$. Therefore,
$f(j,r,n;m)=q^{(r-j)(m-j)}\left[n-m\atop r-j\right]_q.$

Since
\begin{eqnarray*}
f(j,r,n;m)-f(j,r,n;m+1)&=&q^{(r-j)(m-j)}\left[n-m\atop r-j\right]_q-q^{(r-j)(m+1-j)}\left[n-m-1\atop r-j\right]_q\\
&=& (q^{r-j}-1)\frac{q^{(r-j)(m-j)}\prod_{l=n-m-(r-j)+1}^{n-m-1}(q^l-1)}{\prod_{l=1}^{r-j}(q^l-1)}\\
&\geq&0,
\end{eqnarray*}
the desired result follows.
\qed \end{pf}

\begin{thm}\label{thm1.4}
Let $i, d,k,n$ be positive integers  with $\lfloor(d+1)/2\rfloor\leq
i\leq d<k$ and $n-k-\bar{s}(k+d-2i)\geq d-i$. If $k-i\geq2$ and
$1\leq \bar{s}\leq q(q^{k-1}-1)\big/(q^{k-i}-1)$, then the following
hold:

\begin{itemize}
\item[\rm(i)] $M_q(i;d,k,n)$ is an $\bar{s}^{\bar{e}_2}$-disjunct matrix, where
$$\bar{e}_2=q^{(d-i)(k+\bar{s}(k+d-2i)-i)}\left[n-k-\bar{s}(k+d-2i)\atop d-i\right]_q
\left(q^{k-i}\left[k-1\atop
i-1\right]_q-(\bar{s}-1)q^{k-i-1}\left[k-2\atop
i-1\right]_q\right)-1;$$

\item[\rm(ii)] For a given $k$, if $i<d$, then $\lim\limits_{n\longrightarrow +\infty}
\frac{\bar{e}_2+1}{\bar{e}_1+1}=+\infty.$
\end{itemize}
\end{thm}

\begin{pf} (i)
Let $B_{0},B_{1},\ldots,B_{\bar{s}}\in\left[[n]\atop k\right]_q$ be
any $\bar{s}+1$ distinct columns of $M_q(i;d,k,n)$. Clearly, $B_0$
contains $\left[k\atop i\right]_q$ many $i$-dimensional subspaces.
To obtain the maximum number of $i$-dimensional subspaces of $B_0$
in
$$
B_0\cap \bigcup\limits_{j=1}^{\bar{s}}B_j=\bigcup\limits_{j=1}^{\bar{s}}(B_0\cap B_j),
$$
we may assume that $\dim(B_0\cap B_j)=k-1$ for each
$j\in\{1,2,\ldots,\bar{s}\}$. Then each $B_j$ contains
$\left[k-1\atop i\right]_q$ many $i$-dimensional subspaces of $B_0$.
However, any two distinct $B_j$ and $B_l$ intersect at a
$(k-2)$-dimensional subspace. Therefore, only $B_1$ contains
$\left[k-1\atop i\right]_q$ many $i$-dimensional subspaces of $B_0$,
while each of $B_2,B_3,\ldots,B_{\bar{s}}$ contains at most
$\left[k-1\atop i\right]_q-\left[k-2\atop i\right]_q$ many
$i$-dimensional subspaces of $B_0$ not contained in  $B_1$.
Consequently, the number of $i$-dimensional subspaces of $B_0$ not
contained in  $B_1,B_2,\ldots,B_{\bar{s}}$ is at least
\begin{eqnarray*}
\alpha&=&\left[k\atop i\right]_q-\left[k-1\atop i\right]_q-
(\bar{s}-1)\left(\left[k-1\atop i\right]_q-\left[k-2\atop i\right]_q\right)\\
&=&q^{k-i}\left[k-1\atop i-1\right]_q-(\bar{s}-1)q^{k-i-1}\left[k-2\atop i-1\right]_q.
\end{eqnarray*}

Let $D\in\left[[n]\atop d\right]_q$ satisfying  $\dim(D\cap B_0)=i$.
If there exists $j\in \{1,2,\ldots,\bar{s}\}$ such that  $\dim(D\cap
B_j)=i$,   by $(D\cap B_0)+(D\cap B_j)\subseteq D$, we have
\begin{eqnarray*}
&&\dim(B_0\cap B_j)\geq\dim(D\cap B_0\cap B_j)\\
&=&\dim(D\cap B_0)+\dim(D\cap B_j)-\dim((D\cap B_0)+(D\cap B_j))\\
&\geq& 2i-d.
\end{eqnarray*}
Suppose $\dim(B_0\cap B_j)\geq2i-d$ for each $j\in\{1,2,\ldots,\bar{s}\}$.
Then
\begin{eqnarray*}
&&\dim(B_0+B_1+\cdots+B_{\bar{s}})\\
&=&\dim(B_0+B_1+\cdots+B_{\bar{s}-1})+\dim B_{\bar{s}}
-\dim((B_0+B_1+\cdots+B_{\bar{s}-1})\cap B_{\bar{s}})\\
 &\leq&\dim(B_0+B_1+\cdots+B_{\bar{s}-1})+\dim B_{\bar{s}}
-\dim(B_0\cap B_{\bar{s}})\\
&\leq&\dim(B_0+B_1+\cdots+B_{\bar{s}-1})+k+d-2i\\
&\leq&\dim B_0+\bar{s}(k+d-2i)\\
&=&k+\bar{s}(k+d-2i).
\end{eqnarray*}
Let $P$ be a given $i$-dimensional subspace   of $B_0$ not contained
in $B_1,B_2,\ldots,B_{\bar{s}}$. By Lemma~\ref{lem2.2}, the number
of $d$-dimensional subspaces $D$ in $\mathbb F_q^n$ satisfying $
D\cap (B_{0}+B_1+\cdots+B_{\bar{s}}) =P$ is at least
 $$
q^{(d-i)(k+\bar{s}(k+d-2i)-i)}\left[n-k-\bar{s}(k+d-2i)\atop
d-i\right]_q.
$$
Clearly, $D\cap B_0=P$ and $\dim(D\cap B_j)\not=i$ for each
$j\in\{1,2,\ldots,\bar{s}\}$. Therefore, the number of
$d$-dimensional subspaces $D$ in $\mathbb F_q^n$ satisfying
$\dim(D\cap B_{0})=i$ and $\dim(D\cap B_{j})\not=i$ for each
$j\in\{1,2,\ldots,\bar{s}\}$ is at least $$\alpha
q^{(d-i)(k+\bar{s}(k+d-2i)-i)}\left[n-k-\bar{s}(k+d-2i)\atop
d-i\right]_q.$$ Since $\bar{e}_2\geq0$, $\alpha>0$, which implies
that
$$\bar{s}\leq\frac{q^{k-i}\left[k-1\atop i-1\right]_q}{q^{k-i-1}\left[k-2\atop i-1\right]_q}
=\frac{q(q^{k-1}-1)}{q^{k-i}-1}.$$ Hence, (i) holds.

(ii) is straightforward by (i)  and Theorem~\ref{thm1.3}.\qed
\end{pf}

\begin{thm}\label{thm1.5}
Let $i,d,k,n$ be positive integers with $1\leq
i<\lfloor(d+1)/2\rfloor$ and $ d<k,n-(\bar{s}+1)k\geq d-i$. If
$1\leq \bar{s}\leq q(q^{k-1}-1)\big/(q^{k-i}-1)$, then the following
hold:

\begin{itemize}
\item[\rm(i)]
$M_q(i;d,k,n)$ is an $\bar{s}^{\bar{e}_2}$-disjunct matrix,
where
$$\bar{e}_2=q^{(d-i)((\bar{s}+1)k-i)}\left[n-(\bar{s}+1)k)\atop d-i\right]_q
\left(q^{k-i}\left[k-1\atop i-1\right]_q-(\bar{s}-1)q^{k-i-1}\left[k-2\atop i-1\right]_q\right)-1;$$

\item[\rm(ii)] For a given $k$,  $\lim\limits_{n\longrightarrow +\infty}
\frac{\bar{e}_2+1}{\bar{e}_1+1}=+\infty.$
\end{itemize}
\end{thm}

\begin{pf}
The proof   is similar to that of Theorem~\ref{thm1.4}, and will be
omitted. \qed \end{pf}

For $q=2,k=8,n=60$, Table 1 shows the disjunct property of our
designs and Ngo and Du's designs for small $i, d, \bar s.$

  \begin{table} \caption{Disjunct property of $M_q(d,k,n)$ and
$M_q(i;d,k,n)$}
 \begin{center}
\begin{tabular}{|l|l|l|l|l|}
\hline $(i,d)$ & $\bar{s}$ & $\bar{e}_1$&$\bar{e}_2$&Remarks\\
\hline (1,2)&2&6111&$36893488146882232319$& Theorem~ \ref{thm1.4}\\
 \hline (1,3) &2&74927&$3544607988605033156167647492927651839$& Theorem~ \ref{thm1.5}\\
  \hline (2,3) &4&54095&$599519146661432524799$& Theorem~ \ref{thm1.4}\\
 \hline (1,4)&2&177815&$284599986330728289752034695103377217756856319$& Theorem~ \ref{thm1.5}\\
 \hline (2,4)&4&155495&$28799857511436549196854689617936383999$& Theorem~ \ref{thm1.4}\\
 \hline (3,4)&8&110855&$800925501358079$& Theorem~ \ref{thm1.4}\\
   \hline
\end{tabular}
 \end{center}
\end{table}

\small
\section{Concluding remarks}
\begin{itemize}
\item[\rm(i)]
For  given positive integers $d<k$,  $\lim\limits_{n\longrightarrow
+\infty} \frac{\left[n\atop d\right]_q}{\left[n\atop k\right]_q}=0$.
This shows that the test-to-item of $M_q(i;d,k,n)$ is small enough
when $n$ is large enough. By Theorems~\ref{thm1.4} and \ref{thm1.5},
our pooling designs are much better than Ngo and Du's designs when
$n$ is large enough.

\item[\rm(ii)] Ngo \cite{Ngo2} improved the error-tolerant property
of $M_q(d,k,n)$ for $\bar s\geq q+2$, $\bar s\geq q+3$  and $\bar
s\geq q+4$, respectively. By a similar method, we also can improve
the error-tolerant property of $M_q(i; d,k,n)$ for these cases.

\item[\rm(iii)] For positive integers $1\leq d<k< n$, how
about the error-tolerant property of $M_q(0;d,k,n)$?
\end{itemize}

\section*{Acknowledgment}
This research is partially supported by    NSF of China (10971052,
10871027),   NCET-08-0052, Langfang Teachers' College (LSZB201005),
and   the Fundamental Research Funds for
the Central Universities of China.


\begin{thebibliography}{00}\frenchspacing

\bibitem{CD}
Y. Cheng and D. Du, Efficient constructions of disjunct matrices
with applications to DNA library screening, J. Comput. Biol. 14 (2007) 1208-1216.

\bibitem{CD2}
Y. Cheng and D. Du, New constructions of one- and two-stage pooling sesigns,
J. Comput. Biol. 15 (2008) 195-205.

\bibitem{DH}
D. Du and F. K. Hwang, Pooling designs and nonadaptive group testing:
Important Tools for DNA Sequencing, World Scientific, 2006.

\bibitem{DF}
A. G. D'yachkov, F. K. Hwang, A. J. Macula, P. A. Vilenkin and C. Weng,
A construction of pooling designs with some happy surprises, J. Comput. Biol. 12 (2005) 1127-1134.

\bibitem{DF2}
A. G. D'yachkov, A. J. Macula and P. A. Vilenkin,
Nonadaptive and trivial two-stage group testing with
error-correcting $d^e$-disjunct inclusion matrices, In: Entropy, Search, Complexity,
Bolyai society mathematical studied, vol. 16, Spring, Berlin, pp 71-83, 2007.

\bibitem{GW}
J. Guo and K. Wang, A construction of pooling designs with high
degree of error correction, Preprint.


\bibitem{HW}
T. Huang and C. Weng, Pooling spaces and non-adaptive pooling designs,
Discrete Math. 282 (2004) 163-169.

\bibitem{Macula}
A. J. Macula, A simple construction of $d$-disjunct matrices with
certain constant weights, Discrete Math. 162 (1996) 311-312.

\bibitem{Macula2}
A. J. Macula, Error-correcting non-adaptive group testing with
 $d^e$-disjunct matrices, Discrete Appl. Math. 80 (1997) 217-222.

\bibitem{NG}
J, Nan and J. Guo, New error-correcting pooling designs associated
with finite vector spaces, J. Comb. Optim. 20 (2010) 96-100.


\bibitem{Ngo}
H. Ngo and D. Du, New constructions of non-adaptive and
error-tolerance pooling designs, Discrete Math. 243 (2002)
161-170.

\bibitem{Ngo2}
H. Ngo, On a hyperplane arrangement problem and tighter analysis of
an error-tolerant pooling design, J. Comb. Optim. 15 (2008) 61-76.


\bibitem{wanbook}
Z. Wan, Geometry of Classical Groups over Finite
Fields, 2nd edition, Science Press, Beijing/New York, 2002.






\end{thebibliography}
\end{document}